\documentclass[10pt]{article}
\usepackage{amsmath, amssymb, braket, amsthm, graphicx, titlefoot, color}
\usepackage{setspace}
\usepackage[square,numbers]{natbib}
\newtheorem{definition}{Definition}

\title{Control System on Boolean Networks through Semi Tensor Product Approach and Boolean Calculus}
\author{Souma Mazumdar  \\ Department of Theoretical Sciences \\ S. N. Bose National Centre for Basic Sciences \\ Block - JD, Sector - III, Salt Lake City, Kolkata - 700 106 \\ Email: souma.mazumdar@bose.res.in,  Phone: 09903144810}
\date{}

\begin{document}
	\maketitle 
	
	\begin{abstract}
       We have considered a Boolean control network where the state evolution equations depend on past states, controls and first derivatives of a function with respect to controls. Total approach has been the efficient use of matrix semi tensor product and logical operators and logical equations. We have obtained a absolutely new result in considering the derivative terms to play a role to influence the states. We have discussed about  controllability with an example and also optimal control and its significance in control theory.
       
  %     \noindent{\bf Word Count 5600}
	\end{abstract}

	\section{Introduction}
    Boolean Networks were first proposed inspired by the ideas of systems biology which was developed to mimic biological networks at the very cellular level. It was first introduced by Kauffman\cite{kauffman1969metabolic,alberch1994kauffman}, has become powerful tool for describing, analyzing, and simulating cellular networks. Today this is an active field of research drawing attention not from researchers from biology but also from Physics, Control and Systems Sciences. With strong mathematical foundation it has now curved its own niche. With particular focus on its application in control and systems science we try to show how efficient a tool it is in application in control theory by its own Boolean language still making contacts with basic theories and postulates of classical control theory. Almost all the basic motivations of control theory like controllability, observability and optimal control are still derivable in the Boolean regime using the language of Boolean mathematics. It even goes further as it incoprorates the notion of calculus like derivatives and integrals in its own Boolean version.\par
	In this article we primarily discuss about Boolean networks when they are driven by some external parameters called controls. So these networks can be correctly termed as Boolean control networks. In Boolean mathematics every function is a logical function with some logical operators operating between the variables. In control theory we assume a controlled dynamical system which is expressed as a differential equation either in continous or discrete time. In case of Boolean networks the controlled dynamical systems are actually logical dynamical systems parametrized by controls but always in discrete time. As in classical control theory we have a evolution of state variables which are described by some functions which in this case happens to be logical functions. Moreover we can also consider some evolution equations of control variables which are also described by logical functions.\par	
	A wonderful technique exists in this domain where these logical equations can be transformed into algebraic equations. For this we use the logical matrices in case of any logical equation and transition matrices while converting the equations of state evolution. A wonderful formalism has been developed to achieve this by Cheng and his collaborators and a number of interesting papers are written along this line\cite{cheng2014semi,cheng2010analysis,cheng2005semi}. This formalism introduces a new matrix product called the Semi Tensor Product which is a generalisation of conventional matrix product. Using this approach any logical equation can be efficiently converted into its algebraic form which makes mathematical operations possible on them.\par
	As we know in control theory there is a state evolution law which in this case is given by a logical function of state and control variables. But as Boolean calculus is also a developed subject we try to add that feature to our existing state and control evolution laws. In our article we make a novel approach to think of the function describing the state evolution to be a logical function of state varibales, control variables and also first derivatives of another logical function with respect to control variables. We like to see how the network behaves under this new evolution law. Adopting the technique shown by Cheng and his collaborators we convert each of the state evolution equations into algebraic equations analyze them in the context of basic motivations of control theory like controllability and optimal control.\par
	Apart from analysing the characteristics of a Boolean control network we present a very relevant example with some numerical results to illustrate controllability. We have planned to develop a flow of discussions following some basic structure. After a brief introductory discussion we present the basic formalism related to Boolean control networks along with the formalism of Semi Tensor Product which we use as the main tool in our article. We discuss about different types of logical operators and how to convert the logical equations into algebraic equations using logical matrices each reserved for a particular logical operation. As we have have used Boolean derivatives we do a brief discussion on them. We derive a new result for state evolution using Boolean derivatives and Semi Tensor Product approach for a general Boolean network with $n$ state variables and $m$ control variables. Then we illustrate controllability by a relevant example with some numerical results. Towards the end of our discussion we make a brief discussion about the use of optimal control in the developed framework. Then we conclude with a brief conclusion.
	\subsection{Important Notations and basic formalism}
	We adopt similar technique and approach as shown in the book\cite{cheng2010analysis}. \\ \\
	1. $M_{m \times n}$ is the set of $m \times n$ real matrices. \\
	2. Col$_{i}(M)$ is the $i$-th column of matrix $M$;Col$(M)$ is the set of columns of $M$. \\
	3. $\mathcal{D}_{k}:={1,2, \dots, k}$ \\
	4. $\delta^{i}_{n}:=$Col$_{i}(I_{n})$ $i.e.$ it is the $i$-th column of the identity matrix. \\
	5. $\Delta:=$Col$(I_{n})$ \\
	6. $M \in \mathcal{M}_{m \times n}$ is called a logical matrix if Col$(M) \subset \Delta_{m}$ and the set of $m \times n$ logical functions is denoted by $\mathcal{L}_{m \times n}$ \\
	7. Assume $L \in \mathcal{L}_{m \times n}$, then
	$$L=[\delta^{i_{1}}_{m} \;  \delta^{i_{2}}_{m} \dots \delta^{i_{n}}_{m}]$$
	and its shorthand form is 
	$$L=\delta_{m}[i_{1} \; i_{2} \dots i_{n}]$$
	8. A $k$ dimensional vector with all entries equal to $1$ is denoted by
	$$1_{k}:=(1 \; \dots 1)^{T}$$
	9. $A \ltimes B$ is the semi tensor product(STP) of two matrices $A$ and $B$. The symbol $\ltimes$ is mostly omitted and we express $$AB:= A \ltimes B$$
	Here we take the opportunity to briefly introduce the definition of STP. 
	\begin{definition}
		Let $A \in \mathcal{M}_{m \times n}$ and $B \in \mathcal{M}_{p \times q}$. Denote by $t:=$lcm$(n,p)$. Then we define the semi-tensor product(STP) of $A$ and $B$ as 
		\begin{equation}
		A \ltimes B:=(A \otimes I_{\frac{t}{n}})(B \otimes I_{\frac{t}{p}}) \in \mathcal{M}_{(\frac{mt}{n}) \times (\frac{qt}{p})} 
		\end{equation}
	\end{definition}
	It is to be noted when $n=p$, $A \ltimes B=AB$. So the STP is a generalisation of conventional matrix product. STP keeps almost all the major properties of the conventional matrix product unchanged. \\
	We discuss some basic properties of STP. \\
	$a.$ Associative Law :
	\begin{equation}
	A \ltimes (B \ltimes C)=(A \ltimes B)\ltimes C 
	\end{equation} 
	$b.$ Distributive Law :
	\begin{equation}
	\begin{split}
	& (A+B) \ltimes C = A \ltimes C + B \ltimes C \\
	& A \ltimes(B+C)= A \ltimes B + A \ltimes C 
	\end{split} 
	\end{equation}
	$c.$ Transpose :
	\begin{equation}
	(A \ltimes B)^{T}=B^{T} \ltimes A^{T} 
	\end{equation}
	$d.$ Inverse:\\
	If $A$ and $B$ are invertible then
	\begin{equation}
	(A \ltimes B)^{-1}=B^{-1} \ltimes A^{-1} 
	\end{equation}
	$e.$ Let $X \in \mathbb{R}^{t}$ be a column vector. Then for matrix $M$
	\begin{equation}
	X \ltimes M =(I_{t} \otimes M)\ltimes X 
	\end{equation}	
	10. Let $f: \mathcal{B}^{n} \rightarrow \mathcal{B}$ be a boolean function expressed as
	\begin{equation}\label{relation}
	y=f(x_{1}, \dots, x_{n}) 
	\end{equation}
	where $\mathcal{B}=\{0,1\}$. Identifying
	\begin{equation}\label{booldef}
	\begin{split}
	1 = \delta^{1}_{2} = [1 \; 0]^{T}, \; 0 = \delta^{2}_{2}= [0 \; 1]^{T} 
	\end{split} 
	\end{equation}
	It is adopted as a convention to omit the notations of semi tensor product where any multiplication is assumed to be semi tensor product unless explicitly mentioned.
	\subsubsection{Logical Function}
	We define three logical operators as conjuction, disjunction and negation represented by $\wedge$, $\vee$, $\neg$ respectively. These operators operate between booelan variables. The respective equations can be transformed into algebraic equations by logical matrices represented by $M_{c},M_{d},M_{n}$ each reserved for the three operations respectively. \\
	\textit{Example}: \\
	\begin{itemize}
		\item $p \wedge q=M_{c}pq$
		\item $p \vee q=M_{d}pq$
		\item $p \neg q=M_{n}pq$
	\end{itemize}
where $p,q$ are boolean variables and $M_{c},M_{d},M_{n}$ are the respective logical matrices. There is predefined boolean structure of these matrices as given below. \\
$M_{c}=\delta_{2}[1 \; 2 \; 2 \;2], M_{d}=\delta_{2}[1 \; 1\; 1\; 2], M_{n}=\delta_{2}[2 \; 1]$ \\
We define the power reducing matrices as follows. \\
$p^2=M_{r}p$ where $M_{r}=\delta_{4}[1 \; 4]$. The power reduction formula for a product of variables is given by\cite{cheng2010analysis}
\begin{equation*}
(p_{1}\dots p_{i} \dots p_{n})^{2}=\prod_{i=1}^{n}I_{2^{i-1}}\otimes \left[(I_{2}\otimes W_{[2,2^{n-i}]})M_{r}\right]p_{1}\dots p_{i} \dots p_{n}
\end{equation*}
\\
There are some prescribed rules to manipulate the matrices and the variables.
\begin{itemize}
	\item The matrices should always be pulled at the front and the variables should be pushed at the rear. This is achieved by the rule $pM=(I_{2} \otimes M)p$.
	\item Using swap matrices the order of variables can be swapped as follows $pqr=W_{[4,2]}qrp$. When only two variables are swapped it is given as follows
	$pq=W_{[2]}qp$.
\end{itemize}
\subsubsection{Boolean Derivative} Now we introduce briefly the formalism of Boolean derivative as developed in Boolean mathematics. Denote by $M_{f}$ the structure matrix corresponding to a Boolean derivative.If $x=\ltimes_{i=1}^{n} x_{i}$ Then we have\cite{cheng2011matrix},
\begin{equation*}
\begin{split}
\frac{\partial f}{\partial x_{i}} & =M_{\partial_{i}f}x=M_{f}x \oplus M_{f}x_{1} \dots \bar{x}_{i} \dots x_{n} \\ &
=M_{f}x_{1} \dots x_{i} \dots x_{n} \oplus M_{f}x_{1} \dots \bar{x}_{i} \dots x_{n} \\ &
=M_{f}x_{1} \dots (x_{i}+\bar{x}_{i}) \dots x_{n} \\ &
=M_{f}x_{1} \dots (1) \dots x_{n} \\ &
=M_{f}x_{1} \dots \hat{x}_{i} \dots x_{n}
\end{split}
\end{equation*}
where the $\hat{x}_{i}$ denotes the absence of the variable in that position upon differentiation.
\section{Boolean Control networks}
As mentioned at the prelude a Boolean control network\cite{cheng2010analysis} is a Boolean network driven by some parameters called controls. The network consists of $n+m$ nodes where the $n$ nodes serve as states the the remaining $m$ nodes serve as controls. The states undergo change in discrete time following some logical equation where the function describing the change depends on states and controls of the previous time. In case of $\mu$ memory network the state is influenced by values of previous states and controls upto $\mu$th previous time. But in our case we consider the state to depend only on one previous time that is $\mu=1$.\\
In our article we have considered an added feature for the state to depend on the first derivative of a logical function $g$ with respect to control variables. Moreover the controls themselves being dynamical variables of time are considered to undergo change which is captured by a control update function which only depends on control variables of previous time. \\
Let the set of states be denoted by $X=\{x_{1} \dots x_{n}\}$ and the set of controls be denoted by $U=\{u_{1} \dots u_{m}\}$, then the logical function $g$ can be denoted as $g=g(U)$.
\subsection{Mathematical Formalism}
The state update laws are governed by the following equations.
\begin{equation}\label{state update}
\begin{split}
&x_{1}(t+1)=f_{1}\left(x_{1}(t) \dots x_{n}(t),u_{1}(t) \dots u_{m}(t),\frac{\partial g(U)}{\partial u_{1}} \dots \frac{\partial g(U)}{\partial u_{m}}\right) \\ & 
\vdots \\ &
x_{n}(t+1)=f_{n}\left(x_{1}(t) \dots x_{n}(t),u_{1}(t) \dots u_{m}(t),\frac{\partial g(U)}{\partial u_{1}} \dots \frac{\partial g(U)}{\partial u_{m}}\right)
\end{split}
\end{equation}
The control update laws are given by the following equations
\begin{equation}\label{control update}
\begin{split}
&u_{1}(t+1)=g_{1}\left(u_{1}(t) \dots u_{m}(t)\right) \\ &
\vdots \\ &
u_{m}(t+1)=g_{m}\left(u_{1}(t) \dots u_{m}(t)\right)
\end{split}
\end{equation}
Now let us try deriving a general expression for state update equation. This is a new result which we have derived.\\
Converting the first of the equations of the set (\ref{state update}) into algebraic form we have
\begin{equation}
\begin{split}
&x_{1}(t+1)=M_{f_{1}}\prod^{i_{1}}x_{i}\prod^{j_{1}}u_{j}\prod^{k_{1}}\frac{\partial g(U)}{\partial u_{k}}
\end{split}
\end{equation}
Here there are $i_{1}$, $x$ variables which are multipled over. They may be some or all different. Similarly for the $u_{j}$s and the derivatives. Now converting the $i$ th equation of the set (\ref{state update}) we have,
\begin{equation}
\begin{split}
&x_{i}(t+1)=M_{f_{i}}\prod^{i_{i}}x_{i}\prod^{j_{i}}u_{j}\prod^{k_{i}}\frac{\partial g(U)}{\partial u_{k}}
\end{split}
\end{equation}
We have $n$ such equations which are all multiplied over to give a single variable $x$ capturing the state update of the whole network. Therefore we have
\begin{equation}
\begin{split}
x(t+1)&=\prod_{i=1}^{n}x_{i}(t+1) \\ &
=\left[M_{f_{1}}\prod^{i_{1}}x_{i}\prod^{j_{1}}u_{j}\prod^{k_{1}}\frac{\partial g(U)}{\partial u_{k}}\right] \dots \left[M_{f_{n}}\prod^{i_{n}}x_{i}\prod^{j_{n}}u_{j}\prod^{k_{n}}\frac{\partial g(U)}{\partial u_{k}}\right]
\end{split}
\end{equation}
Now all the transition matrices which are structure matrices $M_{f}$ can be pulled to the front following the proposition mentioned earlier. Doing that, the above equation takes the form,
\begin{equation}
\begin{split}
x(t+1)=& M_{f_{1}}(I_{2^{i_{1}+j_{1}+k_{1}}}\otimes M_{f_{2}}) \dots (I_{2^{i_{1} \dots i_{i-1}+j_{1}\dots j_{i-1}+k_{1}\dots k_{i-1}}}\otimes M_{f_{i}}) \dots  (I_{2^{i_{1} \dots i_{n-1}+j_{1}\dots j_{n-1}+k_{1}\dots k_{n-1}}}\otimes M_{f_{n}})  \\ &
\left[\prod^{i_{1}}x_{i}\prod^{j_{1}}u_{j}\prod^{k_{1}}\frac{\partial g(U)}{\partial u_{k}}\right] \dots \left[\prod^{i_{n}}x_{i}\prod^{j_{n}}u_{j}\prod^{k_{n}}\frac{\partial g(U)}{\partial u_{k}}\right] \\ &
=L_{1}\left[\prod^{i_{1}}x_{i}\prod^{j_{1}}u_{j}\prod^{k_{1}}\frac{\partial g(U)}{\partial u_{k}}\right] \dots \left[\prod^{i_{n}}x_{i}\prod^{j_{n}}u_{j}\prod^{k_{n}}\frac{\partial g(U)}{\partial u_{k}}\right] 
\end{split}
\end{equation}
where $L_{1}=M_{f_{1}}(I_{2^{i_{1}+j_{1}+k_{1}}}\otimes M_{f_{2}}) \dots (I_{2^{i_{1} \dots i_{i-1}+j_{1}\dots j_{i-1}+k_{1}\dots k_{i-1}}}\otimes M_{f_{i}}) \dots  (I_{2^{i_{1} \dots i_{n-1}+j_{1}\dots j_{n-1}+k_{1}\dots k_{n-1}}}\otimes M_{f_{n}})$
\\ Leaving $L_{1}$ let us try to simplify the remaining of the expression of state update equation. \\ Let us consider the product $\left[\prod^{i_{1}}x_{i}\prod^{j_{1}}u_{j}\prod^{k_{1}}\frac{\partial g(U)}{\partial u_{k}}\right] \dots \left[\prod^{i_{n}}x_{i}\prod^{j_{n}}u_{j}\prod^{k_{n}}\frac{\partial g(U)}{\partial u_{k}}\right]$. In the above product different $x_{i}$ occur at different positions. So it is not possible to derive a general expression for the swap matrices to bring the like $x_{i}$s together. Let us denote the swap matrix by $W$ which brings like $x_{i}$s together. then the above product can be written as,
\begin{equation*}
\begin{split}
& \left[\prod^{i_{1}}x_{i}\prod^{j_{1}}u_{j}\prod^{k_{1}}\frac{\partial g(U)}{\partial u_{k}}\right] \dots \left[\prod^{i_{n}}x_{i}\prod^{j_{n}}u_{j}\prod^{k_{n}}\frac{\partial g(U)}{\partial u_{k}}\right] \\ &
=Wx_{1}^{p_{1}} \dots x_{n}^{p_{n}}u_{1}^{q_{1}}\dots u_{m}^{q_{m}}\left(\frac{\partial g}{\partial u_{1}}\right)^{s_{1}}\dots \left(\frac{\partial g}{\partial u_{m}}\right)^{s_{m}}
\end{split}
\end{equation*}
where the power over the factors is the number of times they occur in the above product. \\
Using power reducing matrices the above product can be written as
\begin{equation}
\begin{split}
& \left[\prod^{i_{1}}x_{i}\prod^{j_{1}}u_{j}\prod^{k_{1}}\frac{\partial g(U)}{\partial u_{k}}\right] \dots \left[\prod^{i_{n}}x_{i}\prod^{j_{n}}u_{j}\prod^{k_{n}}\frac{\partial g(U)}{\partial u_{k}}\right] \\ &
=WM_{r}^{p_{1}-1}x_{1} \dots M_{r}^{p_{n}-1}x_{n}M_{r}^{q_{1}-1}u_{1}\dots M_{r}^{q_{m}-1}u_{m}M_{r}^{s_{1}-1}\left(\frac{\partial g}{\partial u_{1}}\right)\dots M_{r}^{s_{m}-1}\left(\frac{\partial g}{\partial u_{m}}\right)
\end{split}
\end{equation}
Again the structure matrices can be pulled to the front. Then the above product takes the form
\begin{equation}\label{L2}
\begin{split}
& \left[M_{f_{1}}\prod^{i_{1}}x_{i}\prod^{j_{1}}u_{j}\prod^{k_{1}}\frac{\partial g(U)}{\partial u_{k}}\right] \dots \left[M_{f_{n}}\prod^{i_{n}}x_{i}\prod^{j_{n}}u_{j}\prod^{k_{n}}\frac{\partial g(U)}{\partial u_{k}}\right] \\ &
=WM_{r}^{p_{1}-1} \dots \left(I_{2^{i-1}} \otimes M_{r}^{p_{i}-1}\right) \dots \left(I_{2^{n-1}} \otimes M_{r}^{p_{n}-1}\right)\left(I_{2^{n}}\otimes M_{r}^{q_{1}-1}\right) \dots \left(I_{2^{n+j-1}}\otimes M_{r}^{q_{j}-1}\right)  \\ &\dots \left(I_{2^{n+m-1}}\otimes M_{r}^{q_{m}-1}\right)\left(I_{2^{n+m}}\otimes M_{r}^{s_{1}-1}\right)\dots\left(I_{2^{n+m+k-1}}\otimes M_{r}^{s_{k}-1}\right) \dots \left(I_{2^{n+2m-1}}\otimes M_{r}^{s_{m}-1}\right) \\ & x_{1} \dots x_{n}u_{1} \dots u_{m}\frac{\partial g}{\partial u_{1}} \dots \frac{\partial g}{\partial u_{m}} \\ &
=L_{2}x_{1} \dots x_{n}u_{1} \dots u_{m}\frac{\partial g}{\partial u_{1}} \dots \frac{\partial g}{\partial u_{m}}
\end{split}
\end{equation}
where,
\begin{equation*}
\begin{split}
& L_{2}=WM_{r}^{p_{1}-1} \dots \left(I_{2^{i-1}} \otimes M_{r}^{p_{i}-1}\right) \dots \left(I_{2^{n-1}} \otimes M_{r}^{p_{n}-1}\right)\left(I_{2^{n}}\otimes M_{r}^{q_{1}-1}\right) \dots \left(I_{2^{n+j-1}}\otimes M_{r}^{q_{j}-1}\right)  \\ &\dots \left(I_{2^{n+m-1}}\otimes M_{r}^{q_{m}-1}\right)\left(I_{2^{n+m}}\otimes M_{r}^{s_{1}-1}\right)\dots\left(I_{2^{n+m+k-1}}\otimes M_{r}^{s_{k}-1}\right) \dots \left(I_{2^{n+2m-1}}\otimes M_{r}^{s_{m}-1}\right)
\end{split}
\end{equation*}
Now consider the variables at the end in the equation (\ref{L2}). We have,
\begin{equation}
\begin{split}
&  x_{1} \dots x_{n}u_{1} \dots u_{m}\frac{\partial g}{\partial u_{1}} \dots \frac{\partial g}{\partial u_{m}} \\ &
=x_{1} \dots x_{n}u_{1} \dots u_{m}\left[M_{g}\hat{u}_{1} \dots u_{i} \dots u_{m}\right]\dots \left[M_{g}u_{1}\dots \hat{u}_{i} \dots u_{m}\right] \dots \left[M_{g}u_{1} \dots u_{i} \dots \hat{u}_{m}\right]
\end{split}
\end{equation}
where $M_{g}$ is the structure matrix of the logical function $g$ and we have used the formula for derivative for each of the derivatives. \\
Now again the structure matrices could be pulled to the front. Then we have the above product as,
\begin{equation}\label{L3}
\begin{split}
&  x_{1} \dots x_{n}u_{1} \dots u_{m}\frac{\partial g}{\partial u_{1}} \dots \frac{\partial g}{\partial u_{m}} \\ &
=(I_{2^{n+m}}\otimes M_{g}) \dots (I_{2^{n+im-(i-1)}}\otimes M_{g}) \dots (I_{2^{n+m^2-(m-1)}}\otimes M_{g}) \\ & x_{1} \dots x_{n}u_{1} \dots u_{m}[\hat{u}_{1} \dots u_{i} \dots u_{m}] \dots [u_{1} \dots \hat{u}_{i} \dots u_{m}] \dots [u_{1} \dots u_{2} \dots \hat{u}_{m}] \\ &
=L_{3}x_{1} \dots x_{n}u_{1} \dots u_{m}[\hat{u}_{1} \dots u_{i} \dots u_{m}] \dots [u_{1} \dots \hat{u}_{i} \dots u_{m}] \dots [u_{1} \dots u_{2} \dots \hat{u}_{m}]
\end{split}
\end{equation}
where $L_{3}=(I_{2^{n+m}}\otimes M_{g}) \dots (I_{2^{n+im-(i-1)}}\otimes M_{g}) \dots (I_{2^{n+m^2-(m-1)}}\otimes M_{g})$ \\
Now let us consider the part apart from $L_{3}$ in the equation (\ref{L3}).
We have,
\begin{equation}
\begin{split}
& x_{1} \dots x_{n}u_{1} \dots u_{m}[\hat{u}_{1} \dots u_{i} \dots u_{m}] \dots [u_{1} \dots \hat{u}_{i} \dots u_{m}] \dots [u_{1} \dots u_{2} \dots \hat{u}_{m}] \\ &
=x_{1} \dots x_{n}W_{[2^{m-1},2]}u_{2}\dots u_{m}[u_{1}\dots u_{i} \dots u_{m}][u_{1}\hat{u}_{2} \dots u_{m}] \dots [u_{1}\dots u_{i}\dots u_{m}]
\end{split}
\end{equation}
where we have filled the vacant position of $u_{1}$ in the first bracket by the variable $u_{1}$ in front of the brackets using swap matrix. Now the swap matrix can be brought in front. So we have the above expression as
\begin{equation}
\begin{split}
\left(I_{2^n}\otimes W_{[2^{m-1},2]}\right)x_{1}\dots x_{n}u_{2} \dots u_{m}[u_{1}\dots u_{i} \dots u_{m}][u_{1}\hat{u}_{2} \dots u_{m}] \dots [u_{1}\dots u_{i} \dots \hat{u}_{m}]
\end{split}
\end{equation}
The same trick could be applied for the variable $u_{2}$ to fill up its vacant position in the second bracket. By doing that we obtain the expression as
\begin{equation}
\begin{split}
(I_{2^{n}}\otimes W_{[2^{m-1},2]})(I_{2^{n}} \otimes W_{[2^{2m-1},2]})x_{1}\dots x_{n}u_{3}\dots u_{m}[u_{1} \dots u_{m}][u_{1} \dots u_{m}][u_{1}u_{2}\hat{u}_{3}\dots u_{m}] \dots[u_{1}\dots \hat{u}_{m}]
\end{split}
\end{equation}
When all the vacant positions are filled up the resultant expression stands as follows
\begin{equation}
\begin{split}
(I_{2^{n}\otimes W_{[2^{m-1},2]}})(I_{2^{n}}\otimes W_{[2^{2m-1},2]})\dots (I_{2^{n}}\otimes W_{[2^{m^2-1},2]})x_{1}\dots x_{n}[u_{1}\dots u_{m}] \dots [u_{1}\dots u_{m}]
\end{split}
\end{equation}
where there are $m$ such brackets. So it can be further simplied as
\begin{equation}
\begin{split}
& (I_{2^{n}\otimes W_{[2^{m-1},2]}})(I_{2^{n}}\otimes W_{[2^{2m-1},2]})\dots (I_{2^{n}}\otimes W_{[2^{m^2-1},2]})x_{1}\dots x_{n}[u_{1}\dots u_{m}] \dots [u_{1}\dots u_{m}] \\ & = \prod_{i=1}^{m}(I_{2^{n}}\otimes W_{[2^{im-1},2]})x_{1} \dots x_{n}(u_{1} \dots u_{m})^{m}
\end{split}
\end{equation}
Using the expression for power reduction for a product of variables we can simplify it further to
\begin{equation}
\begin{split}
& (I_{2^{n}\otimes W_{[2^{m-1},2]}})(I_{2^{n}}\otimes W_{[2^{2m-1},2]})\dots (I_{2^{n}}\otimes W_{[2^{m^2-1},2]})x_{1}\dots x_{n}[u_{1}\dots u_{m}] \dots [u_{1}\dots u_{m}] \\ &
=\prod_{i=1}^{m}(I_{2^{n}}\otimes W_{[2^{im-1},2]})x_{1} \dots x_{n}\left[\prod_{i=1}^{m}I_{2^{i-1}}\otimes \left[I_{2}\otimes W_{[2,2^{m-i}]}M_{r}\right]\right]^{m-1}(u_{1}\dots u_{m}) \\ &
=\prod_{i=1}^{m}(I_{2^{n}}\otimes W_{[2^{im-1},2]})\left[I_{2^{n}}\otimes \left[\prod_{i=1}^{m}I_{2^{i-1}}\otimes \left[I_{2}\otimes W_{[2,2^{m-i}]}M_{r}\right]\right]^{m-1}\right]x_{1}\dots x_{m}u_{1}\dots u_{m} \\ &
=L_{4}x_{1}\dots x_{m}u_{1}\dots u_{m}
\end{split}
\end{equation}
where $L_{4}=\prod_{i=1}^{m}(I_{2^{n}}\otimes W_{[2^{im-1},2]})\left[I_{2^{n}}\otimes \left[\prod_{i=1}^{m}I_{2^{i-1}}\otimes \left[I_{2}\otimes W_{[2,2^{m-i}]}M_{r}\right]\right]^{m-1}\right]$
Now we do the final step of bringing all the $u$ variables together in front of all the $x$ variables. For that consider the expression $x_{1} \dots x_{n}u_{1} \dots u_{m}$. \\
We have,
\begin{equation}
\begin{split}
& x_{1}\dots x_{n}u_{1}\dots u_{m} \\ &
=W_{[2,2^{n}]}u_{1}x_{1}\dots x_{n}u_{2} \dots u_{m} \\ &
=W_{[2,2^{n}]}u_{1}W_{[2,2^{n}]}u_{2}x_{1}\dots x_{n}u_{3}\dots u_{m} \\ &
=W_{[2,2^{n}]}(I_{2}\otimes W_{[2,2^{n}]})u_{1}u_{2}x_{1}\dots x_{n}u_{3} \dots u_{m}
\end{split}
\end{equation}
When all the $u$ variables are brought in front of all the $x$ variables and arranged in order we have the final expression as
\begin{equation}
\begin{split}
& x_{1}\dots x_{n}u_{1}\dots u_{m} \\ &
=W_{[2,2^{n}]}(I_{2}\otimes W_{[2,2^{n}]})\dots (I_{2^{i-1}}\otimes W_{[2,2^{n}]})\dots (I_{2^{m-1}}\otimes W_{[2,2^{n}]})u_{1}\dots u_{m}x_{1}\dots x_{n} \\ &
=W_{[2,2^{n}]}\prod_{i=1}^{m-1}\left[I_{2^{i}}\otimes W_{[2,2^{n}]}\right]u_{1}\dots u_{m}x_{1}\dots x_{n} \\ &
=L_{5}u_{1}\dots u_{m}x_{1}\dots x_{n}
\end{split}
\end{equation}
where $L_{5}=W_{[2,2^{n}]}\prod_{i=1}^{m-1}\left[I_{2^{i}}\otimes W_{[2,2^{n}]}\right]$.
So combining all the factors calculated in the preceding steps we have the final state transition equation as
\begin{equation}
\begin{split}
x(t+1)& =L_{1}\ltimes L_{2}\ltimes L_{3}\ltimes L_{4}\ltimes L_{5}\ltimes u_{1}\dots u_{m}x_{1}\dots x_{n} \\ &
=L_{1}\ltimes L_{2}\ltimes L_{3}\ltimes L_{4}\ltimes L_{5}\ltimes ux
\end{split}
\end{equation}
where $L_{1},L_{2},L_{3},L_{4},L_{5}$ are calculated as above and $u=\ltimes_{i=1}^{m}u_{i}$ and $x=\ltimes_{i=1}^{n}x_{i}$.
\section{Controllability}
Controllability is one of the main focus of control theorists while analysing any control problem. Controllability refers to the reachable sets from an initial position with progress of time. Analysing controllability of a control system reveals a great deal about the system such as the trajectory it follows or the fine tuning of the control parameters so that the system can reach a target set in a specified time.\par
In our case we choose a suitable example to illustrate controllability and its importance in control theory. For doing that we have relied on a Matlab package where the computations were done. \\ Consider a system described by the following state and control laws.
\begin{equation}
\begin{split}
& x_{1}(t+1)=x_{2}(t)\vee u_{1}(t)\wedge \frac{\partial g}{\partial u_{2}} \\ &
x_{2}(t+1)=x_{1}(t)\wedge u_{2}(t)\vee \frac{\partial g}{\partial u_{1}} \\ &
u_{1}(t+1)=\neg u_{2}(t) \\ &
u_{2}(t+1)=u_{1}(t) \\ &
g(u_{1},u_{2})=u_{1}\wedge u_{2}
\end{split}
\end{equation}
with the initial condition $x_{1}(0)=x_{2}(0)=u_{1}(0)=u_{2}(0)=\begin{bmatrix}
1 \\ 0
\end{bmatrix}
=\delta_{2}^{1}$
Transforming the logical equations into algebraic equations we have
\begin{equation}
\begin{split}
x_{1}(t+1)& =M_{d}x_{2}(t)u_{1}(t)M_{c}(M_{g}u_{1}u_{2}\oplus M_{g}u_{1}\hat{u}_{2}) \\ &
=M_{d}x_{2}(t)u_{1}(t)M_{c}(M_{g}u_{1}) \\ &
=M_{d}(I_{4}\otimes M_{c})(I_{4}\otimes M_{g})x_{2}(t)u_{1}(t)u_{1}(t) \\ &
=M_{d}(I_{4}\otimes M_{c})(I_{4}\otimes M_{g})x_{2}(t)u_{1}^{2}(t) \\ &
=M_{d}(I_{4}\otimes M_{c})(I_{4}\otimes M_{g})x_{2}(t)M_{r}u_{1}(t) \\ &
=M_{d}(I_{4}\otimes M_{c})(I_{4}\otimes M_{g})(I_{2} \otimes M_{r})x_{2}(t)u_{1}(t)
\end{split}
\end{equation}
where $M_{g}$ stands for the structure matrix of the logical function $g(u_{1},u_{2})$ and the matrices are pulled at the front following the procedure as shown earlier. \\
Converting the second equation we have,
\begin{equation}
\begin{split}
x_{2}(t+1)& =M_{c}(I_{4} \otimes M_{d})(I_{4} \otimes M_{g})x_{1}(t)M_{r}U_{2}(t) \\ &
=M_{c}(I_{4} \otimes M_{d})(I_{4} \otimes M_{g})(I_{2}\otimes M_{r})x_{1}(t)u_{2}(t)
\end{split}
\end{equation}
For the control equations we have,
\begin{equation}
\begin{split}
& u_{1}(t+1)=M_{n}u_{2}(t) \\ &
u_{2}(t+1)=u_{1}(t)
\end{split}
\end{equation}
Now the total state update law of the whole network is given by
\begin{equation}
\begin{split}
x(t+1)&=x_{1}(t+1)x_{2}(t+1) \\ &
=M_{d}(I_{4}\otimes M_{c})(I_{4}\otimes M_{g})(I_{2}\otimes M_{r})x_{2}(t)u_{1}(t)M_{c}(I_{4}\otimes M_{d})(I_{4}\otimes M_{g})(I_{2}\otimes M_{r})x_{1}(t)u_{2}(t) \\ &
=M_{d}(I_{4}\otimes M_{c})(I_{4}\otimes M_{g})(I_{2}\otimes M_{r})\left[I_{4}\otimes[M_{c}(I_{4}\otimes M_{d})(I_{4}\otimes M_{g})(I_{2}\otimes M_{r})]\right]x_{2}(t)u_{1}(t)x_{1}(t)u_{2}(t) \\ &
=M_{d}(I_{4}\otimes M_{c})(I_{4}\otimes M_{g})(I_{2}\otimes M_{r})[I_{4}\otimes[M_{c}(I_{4}\otimes M_{d})(I_{4}\otimes M_{g})(I_{2}\otimes M_{r})]]W_{[2,4]}x_{1}(t)x_{2}(t)u_{1}(t)u_{2}(t) \\ &
=M_{d}(I_{4}\otimes M_{c})(I_{4}\otimes M_{g})(I_{2}\otimes M_{r})[I_{4}\otimes[M_{c}(I_{4}\otimes M_{d})(I_{4}\otimes M_{g})(I_{2}\otimes M_{r})]]W_{[2,4]}x(t)u(t)
\end{split}
\end{equation}
where $x(t)=x_{1}(t)x_{2}(t)$ and $u(t)=u_{1}(t)u_{2}(t)$
Now we solve for $u(t)$ from the initial conditions. We have
\begin{equation}
\begin{split}
u(t+1)&=u_{1}(t+1)u_{2}(t+1) \\ &
=M_{n}u_{2}(t)u_{1}(t) \\ &
=M_{n}W_{[2]}u_{1}(t)u_{2}(t) \\ &
=M_{n}W_{[2]}u(t)
\end{split}
\end{equation}
Starting from the initial conditions we have
\begin{equation}
\begin{split}
& u(1)=M_{n}W_{[2]}u(0) \\ &
u(2)=M_{n}W_{[2]}u(1)=(M_{n}W_{[2]})^{2}u(0) \\ &
\vdots \\ &
u(t)=[M_{n}W_{[2]}]^{t}u(0)
\end{split}
\end{equation}
Therefore,
\begin{equation}
\begin{split}
&x(t+1) =M_{d}(I_{4}\otimes M_{c})(I_{4}\otimes M_{g})(I_{2}\otimes M_{r})[I_{4}\otimes[M_{c}(I_{4}\otimes M_{d})(I_{4}\otimes M_{g})(I_{2}\otimes M_{r})]]W_{[2,4]}x(t)(M_{n}W_{[2]})^{t}u(0)= \\ &
M_{d}(I_{4}\otimes M_{c})(I_{4}\otimes M_{g})(I_{2}\otimes M_{r})[I_{4}\otimes[M_{c}(I_{4}\otimes M_{d})(I_{4}\otimes M_{g})(I_{2}\otimes M_{r})]]W_{[2,4]}[I_{2}\otimes (M_{n}W_{[2]})^{t}]x(t)u(0)= \\ &
M_{d}(I_{4}\otimes M_{c})(I_{4}\otimes M_{g})(I_{2}\otimes M_{r})[I_{4}\otimes[M_{c}(I_{4}\otimes M_{d})(I_{4}\otimes M_{g})(I_{2}\otimes M_{r})]]W_{[2,4]}[I_{2}\otimes (M_{n}W_{[2]})^{t}]W_{[2]}u(0)x(t)
\end{split}
\end{equation}
At this point we consider the initial conditions and also the structure of the structure matrix $M_{g}$. From the function $g(u_{1},u_{2})$ we see that $M_{g}=M_{c}$. Also from the initial conditions $x(0)=x_{1}(0)\ltimes x_{2}(0)=\delta_{2}^{1}\ltimes \delta_{2}^{1}=\delta_{4}^{1}$. Similarly $u(0)=u_{1}(0)\ltimes u_{2}(0)=\delta_{2}^{1}\ltimes \delta_{2}^{1}=\delta_{4}^{1}$. \\
Puting $t=0$ we have,
\begin{equation}
\begin{split}
& x(1)= \\ &
M_{d}(I_{4}\otimes M_{c})(I_{4}\otimes M_{g})(I_{2}\otimes M_{r})[I_{4}\otimes[M_{c}(I_{4}\otimes M_{d})(I_{4}\otimes M_{g})(I_{2}\otimes M_{r})]]W_{[2,4]}[I_{2}\otimes (M_{n}W_{[2]})^{0}]W_{[2]}u(0)x(0)  \\ &
= M_{d}(I_{4}\otimes M_{c})(I_{4}\otimes M_{g})(I_{2}\otimes M_{r})[I_{4}\otimes[M_{c}(I_{4}\otimes M_{d})(I_{4}\otimes M_{g})(I_{2}\otimes M_{r})]]W_{[2,4]}[I_{2}\otimes (M_{n}W_{[2]})^{0}]W_{[2]}\delta_{4}^{1}\ltimes \delta_{4}^{1} \\ &
= L(1)\delta_{4}^{1}\ltimes \delta_{4}^{1}
\end{split}
\end{equation}
where, \\ $L(1)=M_{d}(I_{4}\otimes M_{c})(I_{4}\otimes M_{g})(I_{2}\otimes M_{r})[I_{4}\otimes[M_{c}(I_{4}\otimes M_{d})(I_{4}\otimes M_{g})(I_{2}\otimes M_{r})]]W_{[2,4]}[I_{2}\otimes (M_{n}W_{[2]})^{0}]W_{[2]}$.
A Matlab routine gives
\begin{equation}
\begin{split}
L(1)=&\delta_{4}[1\;1\;1\;2\;2\;2\;2 \;2 \;2 \;2\;2\;2\;2\;2\;2\;2\;2\;2\;2\;2\;2\;2\;2\;2\;2\;2\;2\;2\;2\;2\;2\;2\;1\;1\;1\;2\;2\;2\;2\;2\\ & \;4\;4\;4\;4\;4\;4\;4\;4\;2\;2\;2\;2\;2\;2\;2\;2\;4\;4\;4\;4\;4\;4\;4\;4]
\end{split}
\end{equation}
Also it is found $x(1)=\delta_{4}[1\;1\;1\;2]$. The reachable set is given by columns of $x(1)=\delta_{4}\{1,1,1,2\}$. The distinct sets are $\delta_{4}^{1}$ and $\delta_{4}^{2}$. \\
Similarly puting $t=1$
\begin{equation}
x(2)=L(2)u(0)x(1)
\end{equation}
where,
$L(2)=M_{d}(I_{4}\otimes M_{c})(I_{4}\otimes M_{g})(I_{2}\otimes M_{r})[I_{4}\otimes[M_{c}(I_{4}\otimes M_{d})(I_{4}\otimes M_{g})(I_{2}\otimes M_{r})]]W_{[2,4]}[I_{2}\otimes (M_{n}W_{[2]})^{1}]W_{[2]}$.
As derived from Matlab $L(2)$ comes as,
\begin{equation}
\begin{split}
L(2)=&\delta_{4}[1\;1\;1\;2\;2\;2\;2\;2\;1\;1\;1\;2\;2\;2\;2\;2\;2\;2\;2\;2\;2\;2\;2\;2\;2\;2\;2\;2\;2\;2\;2\;2\;4\;4\;4\;4\;4\;4\;4\;4 \\ & \;2\;2\;2\;2\;2\;2\;2\;2\;4\;4\;4\;4\;4\;4\;4\;4\;2\;2\;2\;2\;2\;2\;2\;2]
\end{split}
\end{equation}
For $x(1)=\delta_{4}^{1}$, $x(2)$ comes as $x(2)=\delta_{4}[1\;1\;1\;2]$. Therefore the reachable set in this case =$\delta_{4}\{1,1,1,2\}$. For $x(1)=\delta_{4}^{2}$, $x(2)$ comes as $x(2)=\delta_{4}[2\;2\;2\;2]$. Therefore the reachable set in this case =$\delta_{4}\{2,2,2,2\}$. Continuing this iteration process we can find the reachable sets of the system for any finite time.
\section{Optimal Control}
Optimal Control refers to a proper control which drives the system into a optimal performance following a optimal path. The focus is on searching for such controls. Here we define a performance index which in this context is termed as the Payoff function. The objective is to optimize the payoff functions by utilizing a proper control. Then the control is called Optimal control. In case of Boolean networks it turns out optimal control drives the system in cycle over time where the system states repeat after a definite period.\\
Consider a Boolean control network with the following state evolution law.
\begin{equation}
x(t+1)=Lu(t)x(t)
\end{equation}
with initial state $x_{0}$.
We define a payoff function as $P(u(t),x(t))$. The average payoff of $x(t,x_{0},u)$ is defined as\cite{cheng2010analysis},
\begin{equation}
J(x(t,x_{0},u))=J(u)=\lim_{t \rightarrow \infty}\frac{1}{T}\sum_{t=1}^{T}P(x(t),u(t))
\end{equation}
The objective of optimal control is to find a control denoted by $u^{\star}$ which maximizes the objective function $J(u)$ that is,
\begin{equation}
J(u^{\star})=\max_{u} J(u)
\end{equation}
It turns out that optimal control drives the system to a loop called a cycle where the system states repeat after definite period. The system as a whole converges to an attractor.
\section{Conclusion}
We like to conclude discussing about the main points of our work and what we have achieved through that. We started our discussion building the main formalism of semi tensor product and the system of equations for a Boolean control network. Boolean networks being logical systems a completely different framework is needed to be build with concepts of logical equations and logical operators which we showed at the prelude. Then considering a general network we obtained a new result on state evolution equations. As mentioned previously it is a totally new approach where we have considered the state update equations to depend not only on states and controls but also derivatives of a function with respect to controls. We successfully implemented our idea to derive the result which is absolutely novel in its structure. Here we like to point that further investigations can be carried out in future to see if any new modifications can be done with the system equations using Boolean calculus. In subsequent flow of our discussion we discussed about controllability, its use in control theory and presented some numerical results. In all our approach we have used extensively the wonderful technique of obtaining matrix products by semi tensor product and also the efficient use of logical operators to reduce the logical equations into algebraic equations. Then towards the end of our discussion we presented a brief notion on optimal control and associated payoff functions and their significance in control theory.
\bibliographystyle{abbrvnat}
\bibliography{library}
\end{document}